\documentclass[12pt]{amsart}
\usepackage{euscript}
\theoremstyle{plain}

\newtheorem{theorem}{Theorem}
\newtheorem{definition}[theorem]{Definition}
\newtheorem{corollary}[theorem]{Corollary}

\newtheorem{example}[theorem]{Example}

\newtheorem{proposition}[theorem]{Proposition}

\def\ev{{\it ev}}
\def\calP{{\EuScript P}}
\def\ext{\mbox{\large$\land$}}
\def\coll#1{\{#1(n)\}_{n\geq 0}}
\def\End{\hbox{${\mathcal E}\hskip -.1em {\it nd}$}}
\def\calI{{\EuScript C}_1}
\def\calJ{{\EuScript J}_1}
\def\Rada#1#2#3{#1_{#2},\dots,#1_{#3}}

\def\Bbb{\mathbb}

\begin{document}
\bibliographystyle{plain}
\pagestyle{plain}

\title{Free Loop Spaces and Cyclohedra}

\author[M. Markl]{Martin Markl}
\thanks{The author was supported by the
        grant GA AV \v CR A1019203
    \newline \indent
        This paper is in final form and no version
        of it will be submitted for publication elsewhere}

\catcode`\@=11
\address{Mathematical Institute of the Academy\\
         \v Zitn\'a 25\\
         115 67 Praha 1\\ Czech~Republic}
\email{markl@math.cas.cz}
\catcode`\@=13


\begin{abstract}
In this note we introduce an action of the module of cyclohedra on the free
loop space. We will then discuss how this action can be used in an
appropriate recognition principle in the same manner as the action of
Stasheff's associahedra can be used to recognize based loop spaces. We will
also interpret one result of R.L.~Cohen as an approximation theorem, in the
spirit of Beck and May, for free loop spaces.
\end{abstract}

\maketitle

Recall that, while a (topological) operad is a
collection $\calP = \coll \calP$ of topological spaces with structure
operations (`compositions')
\begin{equation}
\label{1}
\circ_i : \calP(k) \times \calP(l) \to \calP(k+l-1),\
1\leq i \leq k,
\end{equation}
(see~\cite{may:1972,markl:zebrulka}),
a (right) {\em module} over the operad $\calP$ is a
collection $M = \coll M$ with structure operations (called
again `compositions')
\begin{equation}
\label{2}
\circ_i : M(k) \times \calP(l) \to M(k+l-1),\
1\leq i \leq k,
\end{equation}
satisfying appropriate axioms which are in fact `linearizations' of that of
an operad, see~\cite[Definition~1.3]{markl:zebrulka}.

\begin{example}{\rm\
Let $U$ and $V$ be topological spaces.
The collection
\[
\End_{U,V} := \{{\it Hom}(U^{\times n},V)\}_{n \geq 0}
\]
is a natural right module over the endomorphism operad
$\End_U = \{{\it Hom}(U^{\times n},U)\}_{n \geq 0}$. The structure
operations are given in the expected manner as
\begin{equation}
\label{zase_nestiham} (f \circ_i g)(u_1,\ldots,u_{k+l-1}) := f(\Rada
u1{i-1},g(\Rada ui{l-1+i}),\Rada u{l+i}{k+l-1}),
\end{equation}
for $f: U^{\times k} \to V \in \End_{U,V}(k)$, $g: U^{\times l} \to
U \in \End_U(l)$ and $\Rada u1{k+l-1} \in U$.
}
\end{example}

An algebra over the operad $\calP$ (also called a $\calP$-space) is a
topological space $U$ together with $\Sigma_n$-equivariant maps $A_n:
\calP(n) \to {\it Hom}(U^{\times n},U)$, $n \geq 0$, which are compatible
with the compositions. Thus algebras are `representations' of operads.
Representations of modules are traces, introduced
in~\cite[Definition~2.6]{markl:el} as follows:

\begin{definition}
\label{trace}
Let $M$ be a $\calP$-module as above and $U$ a $\calP$-algebra.
An {\em $M$-trace over $U$} is a topological space $V$
together with $\Sigma_n$-equivariant
maps  $T_n: M(n) \to {\it Hom}(U^{\times n},V)$, $n \geq 0$, which are
compatible with compositions~(\ref{2}) and the $\calP$-algebra
structure on $U$, that is
\[
T_{k+l-1}(m \circ_i p) = T_{k}(m) \circ_i A_{l}(p),
\]
for $m \in M(k)$, $p \in \calP(l)$, $1 \leq i \leq k$, where
$\circ_i$ in the right hand side is as in~(\ref{zase_nestiham}).
\end{definition}

One may also say that a trace is a map $T :M \to \End_{U,V}$ of
collections which is an homomorphism of modules over the homomorphism
$A : \calP \to \End_U$ of operads.

\begin{example}{\rm
Each space $U$ is clearly
a tautological algebra over its endomorphism operad. Similarly, each
space $V$ is a tautological $\End_{U,V}$-trace over
$U$ considered as an $\End_U$-algebra.
}
\end{example}

\begin{example}{\rm
Let $\calI = \coll \calI$ be the little intervals operad (= the little
cubes operad in dimension 1), see~\cite[Example~2.49]{boardman-vogt:73}.
Let $\calJ = \coll \calJ$ be the
collection of `little intervals' in the circle ${\Bbb S}^1$, that is,
$\calJ(n)$ is the space of all linear imbeddings of $n$ unit intervals to
the circle such that the images of interiors are mutually disjoint.
By linear imbeddings we of course mean imbeddings
$\alpha : [0,1] \to {\Bbb S}^1$
that factor as
\[
[0,1] \stackrel{\tilde \alpha}{\longrightarrow} {\Bbb R}
\stackrel{\exp}{\longrightarrow} {\Bbb S^1}
\]
with some linear increasing map $\tilde \alpha$. Therefore, elements of
$\calJ(n)$ are maps
\begin{equation}
\label{win}
d : \bigsqcup_{i=1}^n I_i \to {\Bbb S}^1
\end{equation}
from the disjoint union of $n$-copies of the unit interval $I = [0,1]$ to
the circle such that the restrictions $\alpha_i := d|_{I_i}$ are, for $1
\leq i \leq n$, linear imbeddings and
\[
\alpha_i ({\it int}(I_i)) \cap \alpha_j ({\it int}(I_j)) = \emptyset,
\]
whenever $i \not= j$. Then $\calJ$
is a right $\calI$-module, the module structure being an obvious
generalization of the operadic structure of the little intervals operad.

There is a subspace $F_0({\Bbb S}^1,n) \subset \calJ(n)$ consisting of {\em
cyclically ordered\/} sequences of little intervals, introduced by
R.L.~Cohen in~\cite[Definition~1.2]{cohen:LNM1286}. A map $d$ as
in~(\ref{win}) belongs to $F_0({\Bbb S}^1,n)$ if and only if there is a
sequence
\[
\theta_1 < \theta_2 < \cdots < \theta_n < \theta_1 + 2\pi
\]
of real numbers such that $e^{{\bf i}\theta_i} \in {\it Im}(\alpha_i)$, for
each $1 \leq i \leq n$. The space $F_0({\Bbb S}^1,n)$ carries a natural
right action of the cyclic group ${\Bbb Z}_n$ and $\calJ(n)$ can be
expressed as the induced representation
\begin{equation}
\label{ind}
\calJ(n) \cong {\it Ind}^{\Sigma_n}_{{\Bbb Z}_n} F_0({\Bbb S}^1,n),
\end{equation}
where ${\Bbb Z}_n \subset \Sigma_n$ is interpreted as the subgroup of
cyclic permutations. This equation is  crucial for the proof of
Theorem~\ref{approximation}.
}\end{example}

Let $X$ be a topological space. Recall that the {\em free loop space\/}
$\ext X$ is the space of continuous maps $f : {\Bbb S}^1 \to X$ with
compact-open topology. Given a distinguished point $b \in X$, the {\em based
loop space\/} $\Omega_b X$ or simply $\Omega X$ is the subspace of $\ext X$
of all maps such that $f(*) = b$, where $* := \exp(0)$ is the distinguished
point of ${\Bbb S}^1$.

\begin{theorem}
\label{jsem_v_Tallinnu}
Consider $\Omega X$ as an $\calI$-algebra, with the
classical action of M.~Boardman and
R.~Vogt~\cite[Example~2.49]{boardman-vogt:73}. Then $\ext X$ is a
$\calJ$-trace over the $\calI$-space $\Omega X$.
\end{theorem}

\begin{proof}
A $\calJ$-trace on $\ext X$ is, by
Definition~\ref{trace}, given by assigning,
to any $d \in \calJ(n)$  and based loops $\Rada \phi1n \in
\Omega X$, a free loop $T(d)(\Rada \phi1n) \in \ext X$.

This can be done as follows. Let $d : \bigsqcup_1^n I_n \to {\Bbb S}^1 \in
\calJ(n)$ be as in~(\ref{win}) and denote again by $\alpha_i$ the
restriction $d|_{I_i}$. Then put $T_n(d)(\Rada \phi1n) : {\Bbb S}^1 \to X$
to be $\phi_i \circ \alpha_i^{-1}$ on ${\it Im}(\alpha_i)$, $ 1 \leq i \leq
n$, and $b$ on ${\Bbb S}^1 \setminus {\it Im}(d)$. It is easy to verify that
this indeed defines a trace.
\end{proof}

\begin{example}{\rm
Let $K  = \{K_n\}_{n \geq 0}$ be Stasheff's operad of
associahedra~\cite{stasheff:TAMS63}. It is well-known that $K_n$ can be
interpreted as a compactification of the configuration space of $n$ distinct
points in the unit interval $I$. Let $C_n$ be a similar compactification of
the space of $n$ distinct points in the circle ${\Bbb S}^1$. The space $C_n$
decomposes as $\Sigma_{n-1} \times {\Bbb S}^1 \times W_n$, where
$\Sigma_{n-1}$ is understood here as the quotient of $\Sigma_n$ modulo
cyclic permutations. The space $W_n$ is an $n-1$ dimensional convex polytope
called, by J.~Stasheff~\cite{stasheff:CM97}, the {\em cyclohedron\/}. The
space $W_1$ is the point, $W_2$ is the interval and $W_3$ is the hexagon.
The polyhedron $W_4$ is depicted in Figure~\ref{W4}.
\begin{figure}[tb]
\begin{center}
\unitlength=0.4mm
\begin{picture}(200.11,145.11)

%VIDITELNE HRANY:
\thicklines \put(20.00,10.00){\line(1,0){160.05}}
\put(0.00,20.00){\line(1,1){89.87}} \put(20.00,10.00){\line(-2,1){19.95}}
\put(89.87,109.90){\line(-3,4){9.74}} \put(180.00,10.00){\line(2,1){19.10}}
\put(0.00,20.00){\line(0,1){23.00}} \put(200.00,20.00){\line(0,1){23.00}}
\put(0.00,43.00){\line(1,1){80.00}} \put(100.11,131.11){\line(5,-2){20.00}}
\put(110.00,110.00){\line(3,4){9.89}}
\put(120.00,123.00){\line(1,-1){80.11}}
\put(80.00,123.00){\line(5,2){20.11}}
\put(110.00,110.00){\line(1,-1){90.00}}
\put(110.11,110.00){\line(0,0){0.00}} \put(90.00,110.00){\line(1,0){20.11}}

%NEVIDITELNE HRANY:
\thinlines \put(170.00,25.00){\line(1,2){5.00}}
\put(30.00,25.00){\line(-1,2){5.00}} \put(30.00,25.00){\line(3,1){60.00}}
\put(180.00,10.00){\line(-2,3){10.00}}
\put(85.11,55.11){\line(-3,-1){60.00}} \put(90.00,45.00){\line(1,0){19.89}}
\put(200.00,43.00){\line(-3,-1){15.05}} %\put(175.00,35.00){\line(3,1){6.05}}
\put(109.89,45.11){\line(1,2){4.89}} \put(100.00,65.11){\line(3,-2){15.11}}
\put(90.00,45.11){\line(-1,2){4.89}} \put(85.00,55.00){\line(3,2){15.00}}
\put(20.00,10.00){\line(2,3){10.00}} \put(114.78,55.11){\line(3,-1){60.22}}
\put(0.00,43.00){\line(3,-1){14.96}} %\put(25.00,35.00){\line(-3,1){5.00}}
\put(170.00,25.00){\line(-3,1){60.11}} \put(100.00,65.00){\line(0,1){43}}
\put(100.00,131.00){\line(0,-1){17}}

%KOULE V ROZICH
\put(175.00,35.00){\makebox(0,0)[cc]{$\bullet$}}
\put(200.00,43.00){\makebox(0,0)[cc]{$\bullet$}}
\put(200.00,20.00){\makebox(0,0)[cc]{$\bullet$}}
\put(180.00,10.00){\makebox(0,0)[cc]{$\bullet$}}
\put(170.00,25.00){\makebox(0,0)[cc]{$\bullet$}}
\put(90.00,45.00){\makebox(0,0)[cc]{$\bullet$}}
\put(85.00,55.00){\makebox(0,0)[cc]{$\bullet$}}
\put(100.00,65.00){\makebox(0,0)[cc]{$\bullet$}}
\put(115.00,55.00){\makebox(0,0)[cc]{$\bullet$}}
\put(110.00,45.00){\makebox(0,0)[cc]{$\bullet$}}
\put(80.00,123.00){\makebox(0,0)[cc]{$\bullet$}}
\put(100.00,131.00){\makebox(0,0)[cc]{$\bullet$}}
\put(120.00,123.00){\makebox(0,0)[cc]{$\bullet$}}
\put(90.00,110.00){\makebox(0,0)[cc]{$\bullet$}}
\put(110.00,110.00){\makebox(0,0)[cc]{$\bullet$}}
\put(0.00,43.00){\makebox(0,0)[cc]{$\bullet$}}
\put(0.00,20.00){\makebox(0,0)[cc]{$\bullet$}}
\put(20.00,10.00){\makebox(0,0)[cc]{$\bullet$}}
\put(25.00,35.00){\makebox(0,0)[cc]{$\bullet$}}
\put(30.00,25.00){\makebox(0,0)[cc]{$\bullet$}}

\end{picture}
\end{center}
\caption{The cyclohedron $W_4$.\label{W4}}
\end{figure}
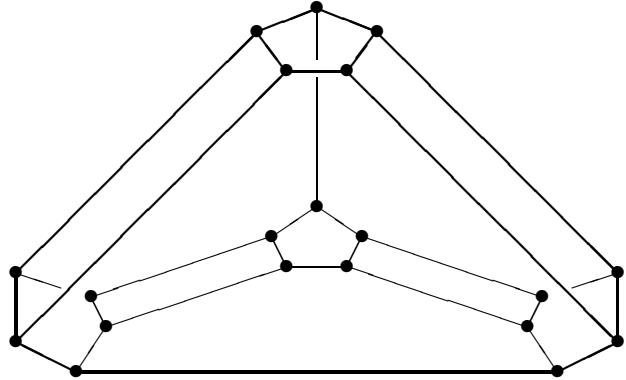
By a slight abuse of terminology we will call $C_n$ the cyclohedron as well.

It is known~\cite[Theorem~2.12]{markl:el} that $C = \{C_n\}_{n \geq 0}$
forms a natural right $K$-module. This structure was analyzed in detail
in~\cite{markl:el}. Philosophically, the couple $(C,K)$ is a minimal
cofibrant  model of $(\calJ,\calI)$ in an appropriate closed model category.
Therefore the following theorem is not surprising. }
\end{example}

\begin{theorem}
\label{jak_dlouho_tady_budeme_bydlet}
Consider $\Omega X$ as a $K$-algebra (= an $A_\infty$-%
space) as in~\cite{stasheff:TAMS63}. Then $\ext X$ is a
$C$-trace over the $K$-algebra $\Omega X$.
\end{theorem}

\begin{proof}
The trace structure can be constructed in the same manner as the operadic
action of the operad $K$ on the based loop space $\Omega
X$~\cite{stasheff:TAMS63} and it is in fact a consequence of the
contractibility of the component of the identity in the space of
endomorphisms of the circle that preserve the base point.

A more conceptual proof follows form the above mentioned fact that the
couple $(C,K)$ is a cofibrant model of $(\calJ,\calI)$ which is more or less
clear from the description of the cyclohedron given in~\cite{markl:el}. Here
`cofibrant' means that the cyclohedron is build up inductively by attaching
cells so the usual obstruction theory for extending maps applies -- we do
not refer to any fancy closed model category structure.

To be precise, let $\rho : K \to \calI$ be a weak homotopy equivalence of
operads. The existence of such a map is folklore and follows from the
cofibrancy of the operad $K$. The map $\rho$ induces a right $K$-module
structure on $\calJ$. Then it follows from the definition of the cyclohedron
and its cellular structure that there exist a weak homotopy equivalence of
right $K$-modules $\lambda : C \to \calJ$.

Let $T : \calJ \to \End_{\Omega X, \land X}$ be the trace structure of
Theorem~\ref{jsem_v_Tallinnu}. Then $T \circ \lambda : C \to \End_{\Omega
X,\land X}$ makes $\ext X$ a $C$-trace over the $K$-algebra $\Omega X$.
\end{proof}

\begin{example}
{\rm\ Let us describe explicitly first components of the structure of
Theorem~\ref{jak_dlouho_tady_budeme_bydlet}. The first piece $T_1 : C_1
\times \Omega X \cong {\Bbb S}^1 \times \Omega X \to \ext X$ of  the trace
is given by the reparametrization
\begin{equation}
\label{3}
T_1(\alpha, \phi)(u) := \phi(u + \alpha),\ \alpha \in {\Bbb S}^1 \cong {\Bbb
R}/{\Bbb N}, \ u \in [0,1].
\end{equation}

To describe a possible choice of the second piece,
$T_2 : C_2 \times \Omega X \times \Omega
X \to \ext X$, introduce coordinates $(\alpha, t)
\in  {\Bbb R}/{\Bbb N} \times I$ of $C_2 \cong {\Bbb S}^1 \times I$.
$T_2$ is then given by
\[
T_2(\alpha, t, \phi, \psi)(u) := (\phi * \psi)(u +
\mbox {$\frac t2$} + \alpha),\ u \in [0,1],
\]
where $\phi * \psi$ is the standard product (loop composition)
in $\Omega X$. The $\Sigma_2$-invariance of $T_2$ follows from the
obvious, but very charming,
equation
\[
(\phi * \psi)(u + \mbox {$\frac 12$}) = (\psi * \phi)(u), \ u \in [0,1],
\]
which holds in $\ext X$ and should be interpreted as a kind of
commutativity for the $*$-product.
}
\end{example}

We believe that an appropriate form of the following  {\em recognition
principle\/} might be true:

\begin{theorem}
\label{zitra_letim_do_Tallinnu}
A couple $(A,B)$ has the homotopy type of $(\ext X, \Omega X)$ if and
only if $B$ is an $A_\infty$-space and $A$ is a $C$-trace
over the $A_\infty$-space $B$.
\end{theorem}

Let $\ext^b X \subset \ext X$ be the subspace of paths $f : {\Bbb S}^1 \to
X$ passing through the distinguished point $b$, that is, $b \in {\it
Im}(f)$. The above recognition principle seems to be related to $\ext^b X$
rather than to $\ext X$. Therefore the following, if fact very surprising,
proposition shared with me by Sa\v sa Voronov, is important.

\begin{proposition}[A.~Voronov]
\label{sasha} Let $X$ be a connected space. Then the inclusion $\ext^c X
\hookrightarrow \ext X$ is, for each $c \in X$, a weak homotopy equivalence.
\end{proposition}

\begin{proof}[Proof due to A.~Voronov]
It is easy to see that the map ${\ev}: \ext^c X \to X$ given by the
evaluation at $* \in {\Bbb S}^1$ is a Hurewicz fibration. The fiber
$\ev^{-1}(b)$ over a point $b \in X$ is the subspace $\Omega_{b}^c X \subset
\Omega_b X$ of loops based at $b$ and passing through $c$. Observe that
$\Omega^c_c X = \Omega_c X$. Consider the following map of fibrations:
\begin{center}
{% Picture saved by xtexcad 2.4
\unitlength=.9pt
\begin{picture}(100.00,105.00)(0.00,-10.00)
\put(100.00,0.00){\makebox(0.00,0.00){$X$}}
\put(0.00,0.00){\makebox(0.00,0.00){$X$}}
\put(100.00,40.00){\makebox(0.00,0.00){$\ext X$}}
\put(0.00,40.00){\makebox(0.00,0.00){$\ext^{c} X$}}
\put(100.00,80.00){\makebox(0.00,0.00){$\Omega_{c} X$}}
\put(0.00,80.00){\makebox(0.00,0.00){$\Omega_c X = \Omega_{c}^c X$}}
\put(100.00,30.00){\vector(0,-1){20.00}}
\put(100.00,70.00){\vector(0,-1){20.00}}
\put(0.00,30.00){\vector(0,-1){20.00}}
\put(0.00,70.00){\vector(0,-1){20.00}} \put(20.00,0.00){\vector(1,0){60.00}}
\put(20.00,40.00){\vector(1,0){60.00}}
\put(40.00,80.00){\vector(1,0){40.00}}
\put(-3.00,20.00){\makebox(0.00,0.00)[r]{\scriptsize $\ev$}}
\put(50.00,45.00){\makebox(0.00,0.00)[b]{\scriptsize \it inclusion}}
\put(103.00,20.00){\makebox(0.00,0.00)[l]{\scriptsize $\ev$}}
\end{picture}}
\end{center}
The proposition follows from the five lemma applied to the induced map of
the long exact homotopy sequences of these fibrations.
\end{proof}

The following corollary, also due to A.~Voronov, immediately follows from
the arguments used in the above proof. This result is so nice and surprizing
that we did not resist the temptation to state it here.

\begin{corollary}[A.~Voronov]
The space $\Omega_{b}^c X$ of based loops passing through a fixed $c \in X$
is homotopy equivalent to the ordinary based loop space $\Omega_b X$.
\end{corollary}

\begin{proof}
The corollary follows from the fact that $\ev: \ext^b X \to X$ is a Hurewicz
fibration, but, once one believes that the statement is true, it is not
difficult write homotopy equivalences explicitly.
\end{proof}

The conceptual meaning of Theorem~\ref{zitra_letim_do_Tallinnu} is the
following. While the recognition principle for based loop spaces (see, for
example,~\cite[Theorem~6.27]{boardman-vogt:73}) says that the path
multiplication, suitably axiomatized, recognizes based loop spaces,
Theorem~\ref{zitra_letim_do_Tallinnu} says that reparametrization~(\ref{3}),
suitably axiomatized, recognizes free loop spaces. The cofibrancy of the
couple $(C,K)$ is very important here, without this assumption one would not
get an `iff' statement.

As argued in~\cite[\S4]{markl-stasheff:ws97}, an {\em approximation
theorem\/} in the spirit of J.~Beck and P.~May~\cite{beck,may:1972} should
describe, for any topological space $X$, the homotopy type of the free
$C$-trace on the $A_\infty$-space $\Omega SX$ which itself has, by the
classical approximation theorem~\cite[Corrolary~6.2]{may:1972}, homotopy
type of the free $A_\infty$-space on $X$. Here, as usual, $S X$ denotes the
reduced suspension.

Let us inspect free traces more closely. Let $\calP$ be an operad and $M$ a
right $\calP$-module as in Definition~\ref{trace}. For a $\calP$-algebra
$A$, the {\em free $M$-trace\/} $T_M(A)$ on $A$ is the quotient of the
disjoint union
\[
\bigsqcup_{n \geq 0} M(n) \times_{\Sigma_n}  A^{\times n}
\]
modulo relations
\begin{eqnarray*}
M(k) \times_{\Sigma_k} A^{\times k}  \ni m \times (\Rada x1{i-1},p(\Rada
xi{i+l-1}),\Rada x{i+l}{k+l-1}) \sim \hskip -2cm&&
\\
&& \hskip -8cm \sim m \circ_i p \times (\Rada x 1{k+l-1}) \in M(k+l-1)
\times_{\Sigma_{k+l-1}} A^{\times k+l-1},
\end{eqnarray*}
for $p \in \calP(l)$, $m \in M(k)$ and $1 \leq i \leq k$. It is easy to see
that $T_M(A)$ is an $M$-trace over the $\calP$-algebra $A$ and that it has
the obvious universal property.

In the special case when $A = \calP X$ is the free $\calP$-algebra on a
space $X$, the above construction can be simplified to a formula which is
very close to~\cite[Construction~2.4]{may:1972}. Namely, $T_M(\calP X)$ is
the quotient of the disjoint union
\[
\bigsqcup_{n \geq 0} M(n) \times_{\Sigma_n}  X^{\times n}
\]
modulo relations
\begin{eqnarray*}
M(n)  \times_{\Sigma_n} X^{\times n} \ni  m \times (\Rada x1{i-1},b,\Rada
x{i+1}n) \sim \hskip -3cm &&
\\
&& \hskip -4cm \sim  m \circ_i *  \times (\Rada x1{i-1},\Rada x{i+1}n) \in
M(n-1)  \times _{\Sigma_{n-1}}  X^{\times n-1},
\end{eqnarray*}
where $b \in X$ and $* \in \calP(0)$ are the distinguished points. Let us
formulate and prove the approximation theorem. Its proof will be based on a
result by R.L.~Cohen.

\begin{theorem}[Approximation theorem for free loop spaces]
\label{approximation} Let $X$ be a connected space. Then there exists the
following diagram of weak homotopy equivalences:
\begin{equation}
\label{casu_je_malo}
{% Picture saved by xtexcad 2.4
\unitlength=1.000000pt
\begin{picture}(140.00,45.00)(0.00,10.00)
\put(80.00,20.00){\makebox(0.00,0.00){\scriptsize $\sim$}}
\put(115.00,50.00){\makebox(0.00,0.00){\scriptsize $\sim$}}
\put(25.00,50.00){\makebox(0.00,0.00){\scriptsize $\sim$}}
\put(70.00,0.00){\makebox(0.00,0.00){$T_{\calJ}(\Omega S X)$}}
\put(140.00,40.00){\makebox(0.00,0.00)[l]{$\ext S X$}}
\put(70.00,40.00){\makebox(0.00,0.00){$T_{\calJ} (\calI X)$}}
\put(0.00,40.00){\makebox(0.00,0.00)[r]{$T_C (KX)$}}
\put(70.00,30.00){\vector(0,-1){20.00}}
\put(100.00,40.00){\vector(1,0){30.00}}
\put(10.00,40.00){\vector(1,0){30.00}}
\end{picture}}
\end{equation}

\vglue6mm
\end{theorem}

\begin{proof}
The vertical homotopy equivalence is induced by the approximation map
$\alpha_1 : \calI X \to \Omega S X$ constructed
in~\cite[Corollary~6.2]{may:1972}.

The left horizontal arrow in~(\ref{casu_je_malo}) is induced by the weak
equivalences $\rho : K \to \calI$ and $\lambda: C \to \calJ$ introduced in
the proof of Theorem~\ref{jak_dlouho_tady_budeme_bydlet}. The existence of
the right horizontal arrow follows from the following identification of
$T_{\calJ} (\calI X)$ with the space $L(X)$ constructed by R.L.~Cohen
in~\cite[Definition~1.3]{cohen:LNM1286}.

It follows from the description of the free trace on the free algebra given
above that $T_{\calJ} (\calI X)$ is the quotient of the disjoint union
\begin{equation}
\label{univ}
\bigsqcup_{n \geq 0} \calJ(n) \times_{\Sigma_n} X^{\times n}
\end{equation}
modulo the relation
\begin{eqnarray*}
\calJ(n) \times_{\Sigma_n} X^{\times n}  \ni (\Rada d1n) \times_{\Sigma_n}
(\Rada x1{i-1},b,\Rada x{i+1}n)\sim \hskip -3cm &&
\\
&& \hskip -9cm \sim (\Rada d1{i-1},\Rada d{i+1}n) \!\times_{\Sigma_n}\!
(\Rada x1{i-1},\Rada x{i+1}n)
 \in  \calJ(n)\times_{\Sigma_{n-1}}\!  X^{\times n-1},
\end{eqnarray*}
It follows from representation~(\ref{ind}) that we may replace~(\ref{univ})
by the disjoint union
\[
\bigsqcup_{n \geq 0} F_0({\Bbb S}^1,n) \times_{{\Bbb Z}_n} X^{\times n}
\]
and modify the relations accordingly. We immediately recognize the space
$L(X)$ of~\cite[Definition~1.3]{cohen:LNM1286}. The right horizontal arrow
of~(\ref{casu_je_malo}) is then identified with the map $h : L(X) \to \ext S
X$ of~\cite[Theorem~1.5]{cohen:LNM1286}. This finishes the proof.
\end{proof}

As proved by Jim Stasheff~\cite{stasheff:TAMS63}, the action of the
associahedra induces an $A_\infty$-structure on the  chain complex
$C_*(\Omega X)$ of the based loop space. Similarly, the action of the
cyclohedra on the free loop space induces a structure of an (algebraic)
strongly homotopy trace on $C_*(\ext X)$. Axioms of this structure are given
more or less explicitly in~\cite{markl:el}.

\vglue2mm \noindent {\bf Acknowledgment.} I am indebted to Jack Morava for
turning my attention to~\cite{cohen:LNM1286}. Theorem~\ref{sasha} and its
corollary is due to Sa\v sa Voronov. I would also like to thank Eugen Paal
for his hospitality during my visit of Tallinn.

%\bibliography{b}

\begin{thebibliography}{1}

\bibitem{beck}
J.~Beck.
\newblock On ${H}$-spaces and infinite loop spaces.
\newblock In {\em Category Theory, Homology Theory and their Applications, III
  (Battelle Institute Conference, Seattle, Wash., 1968, Vol.~3)}, pages
  139--153. Springer, Berlin, 1969.

\bibitem{boardman-vogt:73}
J.M. Boardman and R.M. Vogt.
\newblock {\em Homotopy Invariant Algebraic Structures on Topological Spaces}.
\newblock Springer-Verlag, 1973.

\bibitem{cohen:LNM1286}
R.L. Cohen.
\newblock A model for the free loop space of a suspension.
\newblock {\em Lecture Notes in Mathematics}, 1286:193--207, 1987.

\bibitem{markl:zebrulka}
M.~Markl.
\newblock Models for operads.
\newblock {\em Comm. Algebra}, 24(4):1471--1500, 1996.

\bibitem{markl:el}
M.~Markl.
\newblock Simplex, associahedron, and cyclohedron.
\newblock In J.~McCleary, editor, {\em Higher Homotopy Structures in Topology
  and Mathematical Physics}, volume 227 of {\em Contemporary Math.}, pages
  235--265. Amer. Math. Soc., 1999.

\bibitem{markl-stasheff:ws97}
M.~Markl and J.D. Stasheff.
\newblock On the compactification of configuration spaces.
\newblock {\em Supplem. ai Rend. Circ. Matem. Palermo, Ser. {II}}, 54:83--90,
  1998.

\bibitem{may:1972}
J.P. May.
\newblock {\em The Geometry of Iterated Loop Spaces}, volume 271 of {\em
  Lecture Notes in Mathematics}.
\newblock Springer-Verlag, New York, 1972.

\bibitem{stasheff:TAMS63}
J.D. Stasheff.
\newblock Homotopy associativity of {H-spaces} {I,II}.
\newblock {\em Trans. Amer. Math. Soc.}, 108:275--312, 1963.

\bibitem{stasheff:CM97}
J.D. Stasheff.
\newblock From operads to `physically' inspired theories.
\newblock In J.-L. Loday, J.D. Stasheff, and A.A. Voronov, editors, {\em
  Operads: Proceedings of Renaissance Conferences}, volume 202 of {\em
  Contemporary Math.}, pages 53--81, 1997.

\end{thebibliography}

\def\cprime{$'$}

\end{document}